\documentclass[12pt]{article}

\usepackage{amsmath,amsfonts,amsthm,mathrsfs,amssymb,bbm,leftidx}
\usepackage[all]{xy}
\usepackage{palatino}
\usepackage{a4wide}
\usepackage[
            breaklinks=true,
            pdftitle={An Interpretation of Some Hitchin Hamiltonians In Terms of Isomonodromic Deformation},
            pdfauthor={Michael Lennox Wong}]{hyperref}
\usepackage{color}
\definecolor{tocolor}{rgb}{.1,.1,.5}
\definecolor{urlcolor}{rgb}{.2,.2,.6}
\definecolor{linkcolor}{rgb}{.1,.4,.6}
\definecolor{citecolor}{rgb}{.6,.3,.1}
\hypersetup{colorlinks=true, urlcolor=urlcolor, linkcolor=linkcolor, citecolor=citecolor}

\newcommand{\C}{\mathbb{C}}

\newcommand{\Z}{\mathbb{Z}}

\mathchardef\mhyphen="2D

\renewcommand{\1}{\mathbbm{1}}
\newcommand{\Res}{\text{Res}}

\newcommand{\supp}{\text{supp}}
\newcommand{\pol}{\text{pol}}

\newcommand{\dbar}{\overline{\partial}}
\newcommand{\arsim}{\xrightarrow{\sim}}
\newcommand{\commsq}[8]{\xymatrix{ #1 \ar[r]^{#5} \ar[d]_{#6} & #2 \ar[d]^{#7} \\ #3 \ar[r]_{#8} & #4 } }

\newcommand{\Spec}{\text{Spec}}

\newcommand{\G}{\mathscr{G}}

\renewcommand{\L}{\mathscr{L}}
\newcommand{\M}{\mathscr{M}}

\newcommand{\UU}{\mathfrak{U}}
\newcommand{\Hy}{\mathbb{H}}

\newcommand{\OO}{\mathscr{O}}
\renewcommand{\P}{\mathbb{P}}
\newcommand{\red}{\text{red}}

\newcommand{\Dol}{\text{Dol}}

\newcommand{\g}{\mathfrak{g}}
\newcommand{\h}{\mathfrak{h}}

\renewcommand{\t}{\mathfrak{t}}

\newcommand{\Ad}{\text{Ad}}
\newcommand{\ad}{\text{ad}}

\newcommand{\W}{\mathfrak{W}}
\newcommand{\T}{\mathscr{T}}

\newcommand{\Aut}{\text{Aut}}

\newcommand{\rk}{\text{rk}}

\newcommand{\Y}{\mathcal{Y}}

\newcommand{\PPP}{\mathscr{P}}
\newcommand{\irr}{\text{irr}}
\newcommand{\fr}{\text{cf}}
\newcommand{\fl}{\text{fl}}
\newcommand{\reg}{\text{reg}}

\newcommand{\Sect}{\text{Sect}}
\newcommand{\A}{\mathscr{A}}
\newcommand{\UUU}{\mathcal{U}}

\newcommand{\CC}{\mathscr{C}}
\newcommand{\XX}{\mathscr{X}}

\numberwithin{equation}{section}
\theoremstyle{definition}

\newtheorem{lem}[equation]{Lemma}
\newtheorem{thm}[equation]{Theorem}
\newtheorem{prop}[equation]{Proposition}
\newtheorem{cor}[equation]{Corollary}

\newtheorem{rmk}[equation]{Remark}

\begin{document}

\title{An Interpretation of Some Hitchin Hamiltonians In Terms of Isomonodromic Deformation}

\author{Michael Lennox Wong}

\maketitle

\abstract{This paper deals with moduli spaces of framed principal bundles with connections with irregular singularities over a compact Riemann surface.  These spaces have been constructed by Boalch by means of an infinite-dimensional symplectic reduction.  It is proved that the symplectic structure induced from the Atiyah--Bott form agrees with the one given in terms of hypercohomology.  The main results of this paper adapt work of Krichever and of Hurtubise to give an interpretation of some Hitchin Hamiltonians as yielding Hamiltonian vector fields on moduli spaces of irregular connections that arise from differences of isomonodromic flows defined in two different ways.  This relies on a realization of open sets in the moduli space of bundles as arising via Hecke modification of a fixed bundle.}

\section*{Introduction}

The study of the isomonodromic deformations of connections on holomorphic bundles over Riemann surfaces has its roots in Hilbert's twenty-first problem, or the Riemann--Hilbert problem, which asks whether one can realize a given representation of the fundamental group of a punctured surface as the monodromy of some meromorphic connection whose poles lie at the punctures.  Since the fundamental group is unchanged as we vary the locations of the punctures on the surface, one may seek to determine the precise constraints on these movements ensure that the resulting monodromy remains the same.  This is the problem of isomonodromic deformation.  For simple poles on $\C\P^1$, the answers lie in Schlesinger's equations \cite{Schlesinger}.

For higher order poles, simply defining the monodromy data is a delicate task (see \cite{JMU, Boalchthesis, BoalchGbundles}).  Over $\C\P^1$, the isomonodromy equations were observed by N.\ Hitchin to define a (complex) Poisson manifold, which is symplectic over a dense open set \cite{SchlesingerEquation}.  The symplectic point of view has been further pursued by P.\ Boalch \cite{Boalchthesis, BoalchQHam}, who has described spaces of irregular connections as infinite-dimensional symplectic reductions in the style of Atiyah and Bott \cite{AB}.  Furthermore, using the theory of quasi-Hamiltonian reduction developed by A.\ Alekseev, A.\ Malkin and E.\ Meinrenken \cite{AMM}, he has shown that the space of monodromy data is endowed with a natural symplectic structure, and that the map taking a generic compatibly framed connection to its monodromy data is symplectic.

I.\ Krichever  considered isomonodromic deformation for vector bundles and connections over Riemann surfaces of arbitrary genus $g$ \cite{Krichever}; in this case, the deformation parameters include the moduli of the punctured Riemann surface and the irregular polar part of the connection.  Using the Tjurin parametrization of the moduli space of vector bundles of rank $n$ and degree $n(g-1) + n$, he also represented the flows as non-autonomous Hamiltonian vector fields on the moduli space in the regular singular case.  Following on this, J.\ Hurtubise extended this to bundles of arbitrary degree \cite{Hurtubise}.

We begin with a review of the definition of the monodromy data for a meromorphic connection with a higher order pole at the origin of the unit disc, that is, the Stokes data associated to the connection.  While this is done in \cite{BoalchGbundles}, certain aspects of the construction are used later so we give a short exposition for the sake of clarity.  Section \ref{connections} describes various moduli spaces of bundles with meromorphic connections over a fixed compact Riemann surface, with poles bounded by a fixed divisor, and their symplectic reductions.  As various authors have observed for Higgs bundles \cite{BiswasRamanan, Markman, Bottacin}, and following Hurtubise's description for vector bundles \cite{Hurtubise}, the relevant deformation spaces are given by the first hypercohomology of an appropriate one-step complex; the associated Poisson and symplectic structures are also given in these terms.  The spaces constructed here will be the fibres of a bundle on which the isomonodromy connection will later be described.

The isomonodromy connection is constructed by showing that a connection is determined by its monodromy data, which lie in a space independent of the holomorphic data of the modulus of the Riemann surface or the isomorphism class of the bundle.  Section \ref{globalmonodromy} begins with a review of how the space of monodromy data is constructed.  We cite Boalch's results on the construction of a symplectic form on this space and the fact that the monodromy map, which associates to a triple consisting of a bundle, connection and a compatible framing its monodromy data, is symplectic.  To make the link between Boalch's construction of the moduli spaces as infinite-dimensional symplectic reductions inheriting the Atiyah--Bott symplectic form \cite[\S4]{BoalchQHam} and the hypercohomology realization of the symplectic form given in the previous section \ref{connections}, we justify why these forms agree.

The main results of the paper are given in the final section, but the story told there relies upon being able to realize large open sets in the moduli space of principal bundles as Hecke modifications of a fixed bundle.  Various cases where this is possible were worked out in \cite{HMWCPB}; Section \ref{HM} gives a brief review of this, providing what is necessary for the subsequent discussion.

Section \ref{ID} begins with a description of isomonodromic deformation as a local splitting of, or Ehresmann connection on, a bundle over the space of deformation parameters consisting of the moduli of complex structures on a genus $g$ surface together with a divisor of poles, as well as the irregular part of a connection at the divisor.  The fibres are the spaces of generic compatibly framed connections with fixed irregular part, constructed in Section \ref{connections}.  The rest of the section describes and proves the primary results of this paper.  The main idea is as follows.  Given a tangent vector to the base of the just described bundle (i.e., a deformation of either the modulus of the punctured surface or the irregular part of the connection), the isomonodromy connection produces a unique lift.  Since we are thinking of the moduli of bundles as arising from modifications of a fixed bundle, the isomonodromic deformation of the fixed bundle gives us a second lift.  The difference between these lifts is therefore tangent to the fibre, thus producing a vector field on a moduli space of connections over a fixed Riemann surface and divisor.  A function on this moduli space is then constructed using invariant polynomials, i.e.\ a Hitchin Hamiltonian, which turns out to be a Hamiltonian for the vector field described above.

A first draft of the material appearing here was written as part of a doctoral thesis under the supervision of Jacques Hurtubise.  I thank him warmly for explaining many of the ideas that appear here.  I am also grateful to Marco Gualtieri for his interest in and several discussions on the subject and to Indranil Biswas for clarification on several points.  I would also like to thank Ronnie Sebastian for some troubleshooting help.  Debt is also owed to the referee who caught several inaccuracies, indicated the substance of Remark \ref{quadraticdifferential} and asked for some clarification in the final section.

\section{Local Monodromy} \label{localmonodromy}

Let $G$ be a semisimple complex algebraic group with Lie algebra $\g$, let $T \subseteq G$ be a maximal torus with Lie algebra $\t$, and let $\Phi$ be the associated root system with $\# \Phi =: 2r$.  Let $\Delta \subseteq \C$ be the unit disc with coordinate $z$ and let $P \to \Delta$ be a principal $G$-bundle, necessarily trivial, and let $\nabla$ be a meromorphic connection in $P$ with a pole only at the origin.  In this section, we will briefly review what one needs to obtain the monodromy data, also often referred to as Stokes data, associated to $\nabla$.  As mentioned in the introduction, the definition requires some care and is done by P.\ Boalch in \cite[\S 2]{BoalchGbundles}.  Since we are unlikely to improve upon his exposition, we will describe only what is necessary for our discussion of moduli spaces and refer the reader there for details of the construction.

We will assume that $\nabla$ has a pole of order $k \geq 2$ at $0$.  A \emph{framing} of $P$ at $0$ is a choice of element $s_0 \in P_0$ in the fibre of $P$ above $0$, which we may think of as a section of $P$ over the single point $0$.  A triple $(P, s_0, \nabla)$ will be referred to as a \emph{framed connection}.  Let $s : \Delta \to P$ be a section for which $s(0) = s_0$.  With respect to this section, $\nabla$ becomes a $\g$-valued meromorphic 1-form, and we may consider the lowest order term in the Laurent series expansion:
\begin{align*}
\frac{A_{-k}}{ z^k }.
\end{align*}
The term $A_{-k} \in \g$ depends only on $s_0$ and not on $s$.

A root $\alpha \in \Phi$ may be thought of $\alpha$ as an element of $\t^*$, so that $\ker \, \alpha \subseteq \t$ will be a hyperplane; recall that the set $\t_\reg$ of regular elements of $\t$ is defined to be the complement of all such hyperplanes:
\begin{align*}
\t_\reg := \t \setminus \bigcup_{\alpha \in \Phi} \ker \, \alpha.
\end{align*}
A framed connection $(P, s_0, \nabla)$ is called \emph{compatibly framed} if $A_{-k} \in \t$; it is called \emph{generic} or \emph{non-resonant} if $A_{-k} \in \t_\reg$.

Suppose now that $(P, s_0, \nabla)$ is a generic compatibly framed connection with leading coefficient $A_{-k} \in \t_\reg$.  Then there is a unique formal transformation (i.e.\ transformation in $G( \C [\![ z ]\!])$, so given by a power series which may not converge) whose leading term is the identity with respect to which the connection form is of the form
\begin{align*}
A^0 := \left( \frac{A_{-k}}{ z^k } + \frac{ A_{-(k-1)} }{ z^{k-1} } + \cdots + \frac{ A_{-2} }{ z^2 } + \frac{ \Lambda }{ z } \right) dz,
\end{align*}
where $A_j, \Lambda \in \t, -k \leq j \leq -2$.  $A^0$ is called the \emph{formal type} of $(P, s_0, \nabla)$; the sum of the non-logarithmic terms, i.e.\ $A^0 - \Lambda/z \, dz$, is called the \emph{irregular type}; and $\Lambda$ is called the \emph{exponent of formal monodromy}.

Two compatibly framed connections $(P, s_0, \nabla), (P', s_0', \nabla')$ are said to be \emph{isomorphic} if there exists an isomorphism of $G$-bundles $\varphi : P \to P'$ such that $\varphi(s_0) = s_0'$ and $\varphi^* \nabla' = \nabla$.  In this case, one is generic if and only if the other is, and if they are generic, then they have the same formal type.

Consider the set $\mathscr{H}(A^0)$ of isomorphism classes of generic compatibly framed connections with a fixed formal type $A^0$.  Let $B^+, B^- \subseteq G$ be opposite Borel subgroups containing $T$ and let $U^+, U^-$ be their unipotent radicals.  Given a generic compatibly framed connection $(P, s_0, \nabla) \in \mathscr{H}(A^0)$, there are overlapping sectors in the unit disc on each of which fundamental solutions for the connection (i.e.\ $G$-valued functions $g$ for which the connection form is given by $dg \, g^{-1}$) exist.  To define the Stokes data, one chooses an initial sector, as well as a branch of the logarithm function to specify an initial solution.  On the overlaps of the sectors, the solutions will differ by a constant element of $G$; these elements, the Stokes multipliers, will lie in $U^+$ and $U^-$ for alternate sectors as we go around the disc.  We obtain a mapping
\begin{align*}
\mathscr{H}(A^0) \to (U^+ \times U^-)^{k-1},
\end{align*}
called the \emph{irregular Riemann--Hilbert map}.

\begin{thm} \cite[Theorem 2.8]{BoalchGbundles}
The irregular Riemann--Hilbert map is a bijection.  In particular, $\mathscr{H}(A^0)$ is isomorphic to an affine space of dimension $\# \Phi (k-1) = 2r (k-1)$.
\end{thm}

\section{Connections} \label{connections}

In this section we will be working over a fixed compact Riemann surface $X$ with a fixed effective divisor $D$ of degree $d$.  We will write
\begin{align*}
D & = \sum_{j=1}^m k_j x_j, & D_\red := \sum_{j=1}^m x_j,
\end{align*}
with the $x_j$ distinct so that $d = \sum_{j=1}^m k_j$, $\deg D_\red = m$.

We will let $G, \g, T, \t, \Phi$ be as in Section \ref{localmonodromy}.  By $G_D$ we will mean the group of $G$-valued functions on $D$ (in the schematic sense), so that $G_D = G(\OO_D)$, that is, $G_D$ is the group of $D$-valued points of $G$.  Similarly, the notation of $\g_D$ will often be used.  We will typically think of elements of $\g_D$ as polynomials in local coordinates at the support of $D$ with coefficients in $\g$.

\subsection{Symplectic and Poisson Structures and Reductions}

We consider pairs $(P, \nabla)$, where $P$ is a holomorphic principal $G$-bundle on $X$ and $\nabla$ is a meromorphic connection in $P$ whose poles are bounded by $D$.  In the case where $D$ is reduced, i.e.\ we are considering logarithmic connections, the relevant moduli space can be constructed as in \cite{NitsureConnections} (see also \cite{Simpson94I}); for arbitrary $D$, i.e.\ allowing for irregular poles, the only known construction of the moduli space appears to be an analytic one by an infinite-dimensional symplectic reduction \cite{Boalchthesis,BoalchQHam}.  We will denote by $\L_{X,G}(\varepsilon, D)$ the moduli space of such pairs whose bundle is of topological type $\varepsilon \in \pi_1(G)$, abbreviating to $\L(\varepsilon, D)$ if $X$ and $G$ are understood.  For a pair $(P, \nabla) \in \L(\varepsilon, D)$, the deformation complex is
\begin{align} \label{lsdefmncplx}
\ad \, P \xrightarrow{ -\nabla \cdot} \ad \, P \otimes K(D).
\end{align}
That is, the space of infinitesimal deformations of $(P, \nabla)$ is given by the first hypercohomology group of this complex (cf.\ \cite[Theorem 2.3]{BiswasRamanan}, \cite[Propositions 3.1.2, 3.1.3]{Bottacin}, \cite[Proposition 7.1]{Markman}):
\begin{align*}
\Hy^1( -\nabla \cdot).
\end{align*}

Since the Killing form on $\g$ is $\Ad$-invariant, it gives a well-defined pairing between sections of $\ad \, P$ and
so it follows that $\leftidx{^t}{(-\nabla \cdot)}{} = \nabla \cdot$.  Hence the dual complex to (\ref{lsdefmncplx}) is
\begin{align*}
\ad \, P(-D) \xrightarrow{\nabla \cdot } \ad \, P \otimes K,
\end{align*}
then the diagram
\begin{align} \label{Poisson}
\vcenter{ \xymatrix{
\ad \, P(-D) \ar[r]^{\nabla \cdot} \ar[d]_\1 & \ad \, P \otimes K \ar[d]^{-\1} \\
\ad \, P \ar[r]_{-\nabla \cdot} & \ad \, P \otimes K(D), } }
\end{align}
the top row being the cotangent complex and the bottom the tangent complex, defines a Poisson structure on $\L_{\varepsilon, D}$ \cite[\S6,7]{Markman}.  The vanishing of the Schouten--Nijenhuis bracket, i.e.\ the Jacobi identity, can be proved as in \cite[\S5]{BottacinPoisson} or \cite[\S6]{Polishchuk}.

We will want to realize the spaces $\L(\varepsilon, D)$ slightly differently.  We will consider triples $(P, s, \nabla)$, where $(P, \nabla) \in \L(\varepsilon, D)$ and $s$ is a level structure of $P$ over $D$, i.e.\ a section of $P$ over $D$ or, equivalently, a trivialization of $P$ over $D$; the space of such triples will be denoted $\PPP(\varepsilon, D)$.  Since the space of infinitesimal deformations of a level structure $(P, s)$ is given by $H^1(X, \ad \, P(-D))$, the deformation complex for $(P, s, \nabla)$ is
\begin{align*}
\ad \, P (-D) \xrightarrow{ - \nabla \cdot } \ad \, P \otimes K(D)
\end{align*}
and its dual complex
\begin{align*}
\ad \, P (-D) \xrightarrow{ \nabla \cdot } \ad \, P \otimes K(D).
\end{align*}
Constructing a diagram as in (\ref{Poisson}), since this time we get an isomorphism of complexes, the resulting Poisson structure is non-degenerate and we obtain a symplectic form on $\PPP(\varepsilon, D)$.  Observe that
\begin{align} \label{dimPD}
\dim \PPP(\varepsilon, D) = 2 \dim G(g-1+d).
\end{align}

The space $\PPP(\varepsilon, D)$ admits a free action of $G_D$ with $g \in G_D$ acting on the level structure by
\begin{align*}
g \cdot (P, s, \nabla) = (P, s \cdot g^{-1}, \nabla),
\end{align*}
and it is clear that we may make the identification
\begin{align*}
\L(\varepsilon, D) = \PPP(\varepsilon, D)/ G_D.
\end{align*}
It follows that
\begin{align} \label{dimLD}
\dim \L(\varepsilon, D) = \dim \PPP(\varepsilon, D) - \dim G_D = \dim G \big( 2(g-1) + d \big).
\end{align}

The reason for introducing the level structures is that the symplectic leaves are then easily identified using symplectic reduction (cf.\ \cite[\S6.2]{Markman}).

\begin{prop} \label{LShamaction}
The $G_D$-action on $\PPP(\varepsilon, D)$ is Hamiltonian with moment map $\mu : \PPP(\varepsilon, D) \to \g_D^*$ given by
\begin{align*}
(P, s, \nabla) \mapsto ( s \nabla)_\pol.
\end{align*}
\end{prop}

Here, $(s \nabla)_\pol$ is the Laurent polynomial of $\g$-valued 1-forms we obtain by trivializing $\nabla$ with respect to the section $s$.  It can be paired with an element of $\g_D$ via the invariant bilinear form and taking residues, and hence yields an element of $\g_D^*$.

Thus, the symplectic leaves of $\L(\varepsilon, D)$, which are the symplectic reductions of $\PPP(\varepsilon, D)$, consist of those pairs $(P, \nabla)$ for which the polar part of $\nabla$ lies in a fixed coadjoint orbit in $\g_D^*$.  Again, without a trivialization, the polar part of $\nabla$ is not well-defined, but its coadjoint orbit is.  If $\gamma \subseteq \g_D^*$ denotes a coadjoint orbit, then we will denote the corresponding symplectic leaf in $\L(\varepsilon, D)$ by $\L(\varepsilon, D)^\gamma$.

\subsection{Irregular Parts}

We now consider the subgroup $H_D$ of $G_D$ consisting of elements whose leading term is the identity, i.e., the kernel of the map $G_D = G(\OO_D) \to G(\OO_{D_\red})$.  (This group is referred to as $B_k$ in \cite[\S2]{Boalchthesis} and as $B_D$ in \cite[\S4]{Hurtubise}, but we use $H_D$ so as not to give the impression that we are referring to a Borel subgroup.)  The Lie algebra $\h_D$ of $H_D$ will then be the kernel of $\g_D = \g(\OO_D) \to \g(\OO_{D_\red})$, so if we think of $\g_D$ as polynomials in the local coordinate with coefficients in $\g$, then $\h_D$ is the subalgebra of polynomials with zero constant term.  Dually, if $\g_D^*$ is realized as Laurent polynomials with coefficients in $\g$, then $\h_D^*$ consists of those whose logarithmic term vanishes.

In parallel with Proposition \ref{LShamaction}, we have the following.

\begin{prop}
The $H_D$-action on $\PPP(\varepsilon, D)$ is Hamiltonian with moment map $\mu : \PPP(\varepsilon, D) \to \h_D^*$ given by
\begin{align*}
(P, s, \nabla) \mapsto (s \nabla)_\irr
\end{align*}
where $(s \nabla)_\irr$ is the irregular component of the polar part of $s \nabla$, i.e.\ it is $(s \nabla)_\pol$ with the logarithmic term omitted.
\end{prop}

Observe that the quotient $\L(\varepsilon, D)/ H_D$ is the set of triples $(P, s, \nabla)$ where $s$ is a trivialization of $P|_{D_\red}$, so that in a neighbourhood of each $x_j \in \supp \, D$, we obtain a framed connection; we will denote this quotient by $\L(\varepsilon, D)_\fr$.  The symplectic reductions $\L(\varepsilon, D)_\fr^\gamma$ arising from this action therefore consist of triples $(P, s, \nabla)$, where $s$ is a level structure over $D_\red$ and for which the irregular polar part of $\nabla$ lies in a fixed coadjoint $H_D$-orbit $\gamma \subseteq \h_D^*$.

\subsection{Further Reductions} \label{isoprep}

Let $\W$ be the Weyl group associated to the root system $\Phi$; it may be realized as $\W = N_G(T)/T$, where $N_G(T)$ is the normalizer of $T$ in $G$.  As in Section \ref{localmonodromy}, $r = \frac{1}{2} \# \Phi$ will be the number of positive roots and $l := \rk \, G = \dim T$ will be the rank of $G$ so that $\dim G = 2r + l$.  Let $T_D := T(D) = T(\OO_D)$ be the group of $T$-valued maps on $D$ and $\t_D := \t(\OO_D)$ its Lie algebra.  Let $\PPP(\varepsilon, D, T) \subseteq \PPP(\varepsilon, D)$ be the subspace consisting of triples $(P, s, \nabla)$ for which $(s \nabla)_\pol$ takes values in $\t$ and hence may be considered as an element of $\t_D^*$ and for which $s_\red := s|_{D_\red}$ is a generic compatible framing as in Section \ref{localmonodromy}.

For a fixed $P$ and $\nabla$, it is not hard to see that any two level structures that give elements in $\PPP(\varepsilon, D, T)$ must differ by an element of $N_G(T)_D$, the group of maps from $D$ into $N_G(T)$.  So we get an $N_G(T)_D$-torsor; indeed, we may think of $\PPP(\varepsilon, D, T)$ as a (left) $N_G(T)_D$-bundle over an open set in $\L(\varepsilon, D)$. Therefore, using (\ref{dimLD}),
\begin{align} \label{dimPDT}
\dim \PPP(\varepsilon, D, T) = 2 \dim G (g-1) + (\dim G + l)d = 2 \big( \dim G (g-1) + (r+l)d \big).
\end{align}

Since $T_D \subseteq N_G(T)_D$, it acts on $\PPP(\varepsilon, D, T)$, and as before, we have the following.

\begin{prop}
The $T_D$-action on $\PPP(\varepsilon, D, T)$ is Hamiltonian with moment map $\mu : \PPP(\varepsilon, D, T) \to \t_D^*$ defined as
\begin{align*}
(P, s, \nabla) \mapsto (s \nabla)_\pol.
\end{align*}
\end{prop}

Let us consider the quotient, which we will denote as $\L(\varepsilon, D, T)$, its symplectic leaves $\L(\varepsilon, D, T)^\eta$, and how they compare to those of $\PPP(\varepsilon, D)$.  Elements of $\L(\varepsilon, D, T)$ are triples $(P, w, \nabla)$, where $w$ is a class of level structure with $(w \nabla)_\pol \in \t_D^*$.  There is an induced map
\begin{align*}
\L(\varepsilon, D, T) = \PPP(\varepsilon, D, T) / T_D \to \PPP(\varepsilon, D) / G_D = \L(\varepsilon, D)
\end{align*}
taking
\begin{align*}
(P, w, \nabla) \mapsto (P, \nabla).
\end{align*}
From this expression, it is clear that the fibres are $\W_D$-torsors.

Since the coadjoint orbits in $\t_D^*$ are singletons, a given symplectic leaf $\L(\varepsilon, D, T)^\eta$ of the quotient $\L(\varepsilon, D, T)$ consists of those triples for which $(w \nabla)_\pol = \eta \in \t_D^*$ is fixed (note that this is independent of the representative of $w$).  The preimage of $\L(\varepsilon, D)^\gamma \subseteq \L(\varepsilon, D)$ consists of those $(P, w, \nabla)$ with $(w \nabla)_\pol \in \gamma \cap \t_D^*$.  But this is the union of $\PPP(\varepsilon, D,T)^\eta$ with $\eta \in \gamma \cap \t_D^*$; this intersection is precisely the $\W_D$-orbit of any one of its elements. Thus, the map
\begin{align*}
\bigcup_{\eta \in \gamma \cap \t_D^*} \PPP(\varepsilon, D, T)^\eta \to \L(\varepsilon, D)^\gamma
\end{align*}
is a covering and so an isomorphism on each $\PPP(\varepsilon, D, T)^\eta$.

Similarly, there is an $(T_D \cap H_D)$-action and we record the following.

\begin{prop} \label{torusirregularhamaction}
The $(T_D \cap H_D)$-action on $\PPP(\varepsilon, D, T)$ is Hamiltonian with moment map $\mu : \PPP(\varepsilon, D, T) \to (\t_D \cap \h_D)^*$ given by
\begin{align*}
(P, s, \nabla) \mapsto (s \nabla)_\irr.
\end{align*}
\end{prop}

We may think of an element $g \in N_G(T)_D \cap H_D$ as an $N_G(T)$-valued function on $D$ that is the identity on $D_\red$.  But this means that the image of $g$ must lie in the identity component of $N_G(T)$, which is precisely $T$.  This justifies the following.

\begin{lem}
If $N_G(T)_D = N_G(T)(\OO_D)$, then
\begin{align*}
N_G(T)_D \cap H_D = T_D \cap H_D.
\end{align*}
\end{lem}

\begin{cor}
If $\gamma \subseteq \h_D^*$ is a coadjoint $H_D$-orbit, then $\gamma \cap (\t_D \cap \h_D)^*$ consists of at most a single point.
\end{cor}

From the lemma it follows that the induced map
\begin{align*}
\L(\varepsilon, D, T)_\fr := \PPP(\varepsilon, D, T)/ (T_D \cap H_D) = \PPP(\varepsilon, D, T)/ (N_G(T)_D \cap H_D) \to \PPP(\varepsilon, D)/ H_D
\end{align*}
is an isomorphism onto its image, and that if $\eta \in (\t_D \cap \h_D)^*$ and $\gamma$ is its $H_D$-orbit in $\h_D^*$, then the symplectic reductions can be identified:
\begin{align*}
\L(\varepsilon, D, T)_\fr^\eta = \L(\varepsilon, D)_\fr^\gamma.
\end{align*}
Mixing notation, we will write $\L(\varepsilon, D)_\fr^\eta$ for these spaces, for $\eta \in (\t_D \cap \h_D)^*$.  Since $\dim (T_D \cap H_D) = l(d - m)$, we have
\begin{align} \label{dimLDeta}
\dim \L(\varepsilon, D)_\fr^\eta = 2 \big( \dim G(g-1) + rd + lm \big).
\end{align}

\begin{rmk} \label{etaelts}
One will observe now that for $\eta \in (\t_D \cap \h_D)^*$, elements of $\L(\varepsilon, D)_\fr^\eta$ are triples $(P, s, \nabla)$, where $s$ is a generic compatible framing for $\nabla$, and $\nabla$ is of fixed irregular type at each point of $\supp \, D$.
\end{rmk}

\section{Global Monodromy} \label{globalmonodromy}

For use in this section and the last, we will define a scheme $L$ by
\begin{align} \label{L}
L := \coprod_{j=1}^{m} \Spec \, \C[z]/(z^{k_j}),
\end{align}
so that $L$ is the disjoint union of the $(k_j-1)$th formal neighbourhoods of the origin in $\C$ for $1 \leq j \leq m$.  For now, we will only use $L$ as a way of denoting $m$ points with fixed multiplicities $k_1, \ldots, k_m$, but it will play more of a role in Section \ref{ID}.

\subsection{The Space of Monodromy Data}

We define the manifold of monodromy data following Boalch \cite[\S3]{Boalchthesis} as follows.  For $1 \leq j \leq m$, we set
\begin{align*}
\widetilde{\CC}_j := G \times (U_+ \times U_-)^{k_j-1} \times \t,
\end{align*}
where, in the case $k_j = 1$, we replace $\t$ by the dense open set (though not Zariski open)
\begin{align*}
\t' := \t \setminus \bigcup_{\alpha \in \Phi} \alpha^{-1}(\Z).
\end{align*}
We see that
\begin{align*}
\dim \, \widetilde{\CC}_j = \dim G + (\dim G - l) (k_j - 1) + l = k_j \dim G - l(k_j - 2).
\end{align*}
We now consider the product
\begin{align*}
G^{2g} \times \widetilde{\CC}_1 \times \cdots \times \widetilde{\CC}_{m}
\end{align*}
and observe that it admits a $G$-action:  in each factor $\widetilde{\CC}_j$, $g \in G$ acts by
\begin{align*}
g \cdot ( g_j, \mathbf{K}^j, \Lambda_j) = ( g_j g^{-1}, \mathbf{K}^j, \Lambda_j ),
\end{align*}
and in each factor of $G^{2g}$, the action is by conjugation.  If $\widetilde{\XX}_{g, L}^0$ refers to the submanifold of the product satisfying
\begin{align} \label{monodromyrelation}
[A_1, B_1] \cdots [A_g, B_g] g_1 \exp( 2\pi i \Lambda_1) g_1^{-1} \cdots g_{m} \exp( 2\pi i \Lambda_{m}) g_{m}^{-1} = e,
\end{align}
then the space of monodromy data is then defined to be
\begin{align*}
\XX_{g,L} :=  G \backslash \widetilde{\XX}_{g,L}^0.
\end{align*}
Its dimension is given by
\begin{align*}
\sum_{j=1}^{m} \big( d_j \dim G - l(k_j - 2) \big) + 2g \dim G - 2\dim G = 2 [ \dim G(g-1) + rd + lm].
\end{align*}
Comparing with (\ref{dimLDeta}), we observe that this is precisely $\dim \L(\varepsilon, D)_\fr^\eta$.

Without delving into the theory of quasi-Hamiltonian $G$-spaces and reduction developed in \cite{AMM}, which provides a variation of the well-known theory of symplectic reduction where the moment maps are $G$-valued, we give a brief and rough explanation of how it gives a more geometric construction of $\XX_{g,L}$ together with a holomorphic symplectic form.  The spaces $G^2 = G \times G$ can be thought of as spaces of representations of the fundamental group of a punctured torus; they are quasi-Hamiltonian $G$-spaces \cite[Proposition 3.2]{AMM}.  Glueing tori together corresponds to what is known as the fusion product \cite[\S6]{AMM}, so the data for representations of a genus $g$ surface will come from a $g$-fold fusion product $G^2 \circledast \cdots \circledast G^2$.  Boalch shows that the spaces $\widetilde{\CC}_i$ are quasi-Hamiltonian $(G \times T)$-spaces \cite[Theorem 5]{BoalchQHam}, so that the fusion product
\begin{align*}
G^2 \circledast \cdots \circledast G^2 \circledast \widetilde{\CC}_1 \circledast \cdots \circledast \widetilde{\CC}_m
\end{align*}
is a quasi-Hamiltonian $(G \times T^m)$-space, whose $G$-reduction is precisely $\XX_{g,L}$ as described above.  This has the following consequence.

\begin{thm} \cite[Theorems 3, 4, 5]{BoalchQHam}
The manifold $\XX_{g,L}$ carries a holomorphic symplectic form.
\end{thm}

\subsection{The Monodromy Map} \label{tentacles}

We return to the spaces $\L(\varepsilon, D)_\fr^\eta$ described at the end of Section \ref{isoprep}.  As pointed out in Remark \ref{etaelts}, an element is a triple $(P, s, \nabla)$, where $P$ is a $G$-bundle, $s$ a trivialization over $D_\red$, and $\nabla$ a connection with poles bounded by $D$ and such that $(s \nabla)_\pol = \eta \in (\t_D \cap \h_D)^*$ is fixed.  These spaces are precisely those constructed by Boalch via an infinite-dimensional symplectic reduction \cite[Definition 15, Theorem 9]{BoalchQHam} (cf.\ \cite[\S4,5]{Boalchthesis}), as such they are endowed with a complex analytic  symplectic form, which we will call the Atiyah--Bott form, as it is induced from a symplectic form on a space of $C^\infty$ connections (cf.\ \cite[\S9]{AB}).

Section \ref{localmonodromy} indicated how to define monodromy data at each pole.  There we saw that the data needed to define the monodromy data of a meromorphic connection in the neighbourhood of a single pole was:
\begin{enumerate}\addtolength{\itemsep}{-0.25\baselineskip} \label{localmonodromychoices}
\item a choice of a coordinate $z$ at the pole;
\item a choice of a generic compatible framing;
\item a choice of initial sector $\Sect_1$; and
\item a choice of branch of $\log z$.
\end{enumerate}
At each $x_j \in \supp \, D$, $s$ gives a generic compatible framing, but we will have to choose a coordinate $z_j$ centred at $x_j$, a branch of $\log z_j$ and an initial sector $\Sect_1^j$.  If $z_j'$ is another such choice of coordinate with
\begin{align*}
z_j' = z_j + O(z_j^{k_j+1})
\end{align*}
then the polar part of the expression of $\nabla$ with respect to $s_{x_j}$ is unchanged, so we need only choose an equivalence class of such coordinates, which Krichever defines to be a $k_j$-jet, denoted $[z_j]_{k_j}$ \cite[\S3]{Krichever}.

We would now like to associate global monodromy data to a triple $(P, s, \nabla)$.  For global data, we would like to record the monodromy over a set of $a_i$- and $b_i$-cycles, $1 \leq i \leq g$, that do not intersect the poles, which we will now fix; we would also like to see how the local solutions near the poles relate to each other.  To make sense of this, we choose what Boalch calls a set of ``tentacles'' \cite[Definition 3.9]{Boalchthesis}.  This involves choices of the following data:
\begin{enumerate}\addtolength{\itemsep}{-0.25\baselineskip} \label{tentacle}
\item a base point $y \in X \setminus D$;
\item a point $y_j$ in the coordinate neighbourhood of $x_j$ so that it does not lie on an anti-Stokes direction; we may label the sectors so that $y_j$ lies in $\Sect^j_1$;
\item a branch of $\log z_j$, which we may continue analytically starting in $\Sect^j_1$;
\item a path $\gamma_j$ joining $y$ to $y_j$.
\end{enumerate}
We do this in such a way that the $\gamma_j$ do not intersect any of the $a_i$- or $b_i$-cycles.  To see how and why this can be done, we can cut $X$ along the $a_i$- and $b_i$-cycles to obtain a $4g$-gon with the poles in the interior.  We may assume that $y$ is the base point, which corresponds to the vertices of the $4g$-gon, and then we can choose non-intersecting $\gamma_j$ from $y$ to each $y_j$.  We remark that the choice of tentacles includes a choice of initial sector near each pole as well as the branch of the logarithm there.  Therefore, given $(P, s, \nabla)$, we require the data of a choice of $k_j$-jet near each $x_j$.

Finally, we define the monodromy data as follows.  Since $y$ is a regular point for the connection, there exists a fundamental solution $\phi_0$ in a neighbourhood of $y$, unique up to a constant element of $G$.  We parallel transport it along $\gamma_j$ to $y_j$; then if $\phi_1^j$ is the canonical fundamental solution in $\Sect_1^j$, we have
\begin{align*}
\phi_0 = \phi_1^j g_j
\end{align*}
for some $g_j \in G$.  This $g_j$ is the data that goes into the $G$-component of $\widetilde{\CC}_j$. The Stokes multipliers $K_1^j, \ldots, K_{2(k_j-1)}^j$ go into the factor $(U_+ \times U_-)^{k_j - 1}$ and the exponent of formal monodromy $\Lambda_j$ goes into the $\t$-component.  If $A_i$ is the monodromy around $a_i$ and $B_i$ that around $b_i$, then the $G^{2g}$ factor is
\begin{align*}
(A_1, \ldots, A_g, B_1, \ldots, B_g).
\end{align*}

We therefore get a map $\L(\varepsilon, D)_\fr^\eta \to \widetilde{\CC}_1 \times \cdots \times \widetilde{\CC}_{m} \times G^{2g}$:
\begin{align*}
(P, s, \nabla) \mapsto \left( \prod_{j=1}^{m} (g_j, K_1^j, \ldots, K_{2k_j-2}^j, \Lambda_j), A_i, B_i \right).
\end{align*}
The image clearly consists of elements satisfying (\ref{monodromyrelation}).  To make things independent of the choice of fundamental solution $\phi_0$ in a neighbourhood of $y$, we quotient out by $G$.  This gives us a well-defined map
\begin{align} \label{monodromymap}
\nu : \L_{\varepsilon, D, \fr}^\eta \to \XX_{g,L}
\end{align}
called the \emph{monodromy map}.

\begin{prop} \label{monodromysymplectic}
The monodromy map is a complex analytic symplectomorphism onto its image.
\end{prop}

The proof follows that of Proposition 3.7 in \cite{Boalchthesis} (cf.\ Lemma 3.2 in \cite{Krichever}).  It is sufficient to see that it is an injective map of complex manifolds of the same dimension.  To check injectivity, if $(P, t, \nabla), (P', t', \nabla')$ have the same Stokes data, with canonical fundamental solutions $\phi_j, \phi_j'$ in the Stokes sectors, then $\phi_j' \phi_j^{-1}$ give isomorphisms on the sectors and agree on overlaps.  Furthermore, it is bounded, so extends to an isomorphism over the poles.  The fact that the symplectic structure is preserved follows from \cite[Theorem 9]{BoalchQHam}.

\subsection{Comparison of Symplectic Forms}

The careful reader will have noticed that while we have described the symplectic form on $\L(\varepsilon, D)_\fr^\eta$ in terms of hypercohomology, Boalch's result stating that the monodromy map is symplectic (Proposition \ref{monodromysymplectic}) uses the Atiyah--Bott symplectic form.  We now justify that these are the same.

Let us fix $(P, s, \nabla) \in \L(\varepsilon, D)_\fr^\eta$ and think about what a deformation should look like.  The deformation of complex of $(P, s, \nabla)$ as an element of $\L(\varepsilon, D)$, i.e.\ we are not fixing the irregular part, is
\begin{align*}
\ad \, P(-D_\red) \xrightarrow{-\nabla \cdot} \ad \, P \otimes K(D).
\end{align*}
In terms of \v{C}ech representatives with respect to an open cover $\UU = \{ U_\alpha \}$, a deformation is an element
\begin{align*}
(\tau, \nu) \in Z^1(\UU, \ad \, P(-D_\red)) \oplus C^0(\UU, \ad \, P \otimes K(D))
\end{align*}
satisfying the hypercohomology cocycle condition:
\begin{align} \label{hypercocycle}
- \nabla \tau_{\alpha \beta} = \nu_\alpha - \nu_\beta.
\end{align}
To stay within the symplectic leaf $\L(\varepsilon, D)_\fr^\eta$, we wish to keep constant the irregular part of $\nabla$ and we wish to vary the logarithmic term of $\nabla$ with respect to the framing $s$ by elements of $\t$.  Therefore we wish to restrict $\nu$ to lie in the subsheaf of $\ad \, P(D_\red)$ such that the logarithmic term with respect to $s$ lies in $\t$.  Let us denote this sheaf by $\ad \, P \otimes K(D_\red)^{s,\t} \subseteq \ad \, P \otimes K(D_\red)$.  The condition (\ref{hypercocycle}) means that $\tau$ must lie in the preimage of $\ad \, P \otimes K(D_\red)^{s,\t}$ under $-\nabla \cdot$; let us denote this preimage by $\ad \, P \otimes K(-D_\red)^{s,\t}$.  Therefore the deformation complex should be
\begin{align} \label{deformationcomplexleaf}
\ad \, P(-D_\red)^{s,\t} \xrightarrow{-\nabla \cdot} \ad \, P \otimes K(D_\red)^{s,\t}.
\end{align}

Since we want to compare the hypercohomology form to one given by a symplectic reduction of $C^\infty$ objects, it makes more sense to consider hypercohomology representatives in terms of a Dolbeault resolution.  Since $\dbar$ and $\nabla$ anti-commute, the relevant double complex is
\begin{align} \label{Dolbeaultresolution}
\xymatrixcolsep{4pc}
\vcenter{ \commsq{ C^\infty \big( \ad \, P(-D_\red)^{s,\t} \big) }{ C^\infty \big( \ad \, P \otimes K(D_\red)^{s,\t} \big) }{ \Omega^{0,1} \big( \ad \, P(-D_\red)^{s,\t} \big) }{ \Omega^{0,1} \big( \ad \, P \otimes K(D_\red)^{s,\t} \big) }{ - ( \partial + [A, \cdot ]) }{ - \dbar }{ \dbar }{ - ( \partial + [A, \cdot ]) = - \nabla \cdot } }
\end{align}
The cocycle condition for
\begin{align*}
(\tau, \nu) \in \Omega^{0,1} \big( \ad \, P (-D_\red)^{s,\t} \big) \oplus C^\infty \big( \ad \, P \otimes K(D)^{s,\t} \big),
\end{align*}
is then
\begin{align} \label{Dolbeaultcocycle}
-\nabla \tau = \dbar \nu.
\end{align}
In these terms, the symplectic form defined as in (\ref{Poisson}) will be given explicitly by
\begin{align*}
\omega_\Hy \big( ( \tau_1, \nu_1) , (\tau_2, \nu_2) \big) = \int_X \kappa( \tau_1 , \nu_2) - \kappa( \tau_2 , \nu_1),
\end{align*}
where $\kappa$ is the Killing form on $\g$ appropriately extended to $\ad \, P$-valued forms.

The construction of $\L(\varepsilon, D)_\fr^\eta$ via an infinite-dimensional symplectic reduction can be described as follows.  Consider $P$ as a $C^\infty$ bundle and let $\A(\eta)$ be the set of connections on $P$ whose Taylor expansion with respect to a smooth trivialization extending the compatible framing $s$ is simply $\eta$, thought of as a Laurent polynomial (in fact, a polynomial in $1/z$) with values in $\t \subseteq \g$; let $\A(\eta)_\fl$ denote the subspace of $\A(\eta)$ consisting of flat connections.  If $\nabla \in \A(\eta)_\fl$, then $\nabla^{0,1}$ is non-singular everywhere and so defines a holomorphic structure on $P$.  A holomorphic frame for this holomorphic structure can be obtained by solving $g^{-1} (\dbar g) = \nabla^{0,1}$ (cf.\ \cite[proof of Proposition 4.3]{Boalchthesis}).  With respect to such a frame, the meromorphic connection is given by $\nabla^{1,0}$ (the fact that it is meromorphic follows from flatness).

Fixing an arbitrary $\nabla \in \A(\eta)$, its tangent space is given by
\begin{align*}
T_\nabla \A(\eta) = \left\{ \phi \in \Omega^1 \big( \ad \, P(D) \big) \, \bigg| L_i(\phi) \in \t \frac{dz}{z}, \ 1 \leq i \leq m \right\},
\end{align*}
where $L_i(\phi)$ denotes the Taylor expansion of $\phi$ at the pole $x_i$ with respect to a trivialization in which $\nabla$ is represented by $\eta$ \cite[\S4]{BoalchQHam}.  On this space, the Atiyah--Bott form is given by
\begin{align*}
\omega_{AB}( \phi, \psi) = \int_X \kappa( \phi, \psi).
\end{align*}

The space $\A(\eta)$ is acted upon by the subgroup $\G_1$ of the $C^\infty$ gauge group of $P$ whose Taylor expansion at any pole is the identity.  This action is Hamiltonian with moment map given by the taking of the curvature of a connection.  By definition, the symplectic quotient $\L(\varepsilon, D)_\fr^\eta$ is then the quotient $\A(\eta)_\fl / \G_1$.  Therefore a deformation of $(P, s, \nabla) \in \L(\varepsilon, D)_\fr^\eta$ is represented by an element of $T_\nabla \A(\eta)_\fl$.  The condition for $\phi \in T_\nabla \A(\eta)$ to lie in $T_\nabla(\eta)_\fl$ is
\begin{align*}
0 = F_{\nabla + \epsilon \phi} = d( A + \epsilon \phi) + ( A + \epsilon \phi) \wedge (A + \epsilon \phi) = F_\nabla + \epsilon \nabla \phi,
\end{align*}
where $A$ is a connection form representing $\nabla$.  Thus,
\begin{align*}
T_\nabla \A(\eta)_\fl = \{ \phi \in T_\nabla \A(\eta) \, | \, \nabla \phi = 0 \}.
\end{align*}
Let $\phi = \phi^{1,0} + \phi^{0,1} \in T_\nabla \A(\eta)_\fl$ be the decomposition into $(1,0)$ and $(0,1)$ parts; then choosing a holomorphic frame for $P$ so that $A = A^{1,0}, A^{0,1} = 0$, we obtain
\begin{align*}
0 = d \phi + [A, \phi] = \partial \phi^{0,1} + \dbar \phi^{1,0} + [A, \phi^{0,1}] = \nabla \phi^{0,1} + \dbar \phi^{1,0}.
\end{align*}
But note in this case that $\phi^{1,0}$ is a smooth section of $\ad \, P(D_\red)^{s,\t}$.  This relation shows $\phi^{0,1}$ is a $(0,1)$-form with values in the preimage of $\ad \, P(D_\red)^{s,\t}$, and the relation itself is precisely that of (\ref{Dolbeaultcocycle}).  If $Z^1(-\nabla \cdot)_\Dol$ denotes the space of hypercohomology cocycles for a Dolbeault resolution of the deformation complex, then one has a map $T_\nabla \A(\eta)_\fl \to Z^1(-\nabla \cdot)_\Dol$ simply given by
\begin{align*}
\phi \to (\phi^{0,1}, \phi^{1,0}).
\end{align*}

The Lie algebra of $\G_1$ can be identified as the smooth sections of $\ad \, P$ whose Taylor series at $D$ vanishes.  If $\mu$ is such a section, then the infinitesimal action of $\mu$ at $\nabla$ is readily computed to be $- \nabla \mu \in T_\nabla (\eta)_\fl$.  Using (\ref{Dolbeaultresolution}), it is easy to verify that under the map of the previous paragraph, the subspace of $T_\nabla \A(\eta)_\fl$ generated by the infinitesimal action maps into the space $B^1(-\nabla \cdot)_\Dol$ of coboundaries.  The induced map of the respective quotients is simply the identity map of $T_{(P,s, \nabla)} \L(\varepsilon, D)_\fr^\eta$.

Using these identifications, we can simply check that
\begin{align*}
\omega_{AB}( \phi, \psi) & = \int_X \kappa( \phi^{0,1} + \phi^{1,0}, \psi^{0,1} + \psi^{1,0}) = \int_X \kappa( \phi^{1,0}, \psi^{0,1}) - \kappa( \psi^{1,0}, \phi^{0,1}) \\
& = \omega_\Hy \big( (\phi^{0,1}, \phi^{1,0}), (\psi^{0,1}, \psi^{1,0}) \big)
\end{align*}
noting that $\kappa( \phi^{0,1}, \psi^{0,1})$ and $\kappa( \phi^{1,0}, \psi^{1,0})$ are $(0,2)$- and $(2,0)$-forms, respectively.

\begin{prop}
The Atiyah--Bott and hypercohomology symplectic forms, $\omega_{AB}$ and $\omega_\Hy$, on $\L(\varepsilon, D)_\fr^\eta$ agree.
\end{prop}

\begin{rmk}
The fact that the monodromy map is symplectic indicates that the symplectic structure on the moduli spaces $\L(\varepsilon, D)_\fr^\eta$ is expressible in terms independent of the analytic data of the isomorphism class of the bundle, the connection or even the complex structure of the Riemann surface.  The hypercohomology perspective almost gives another way to see this.  We may consider the subsheaf $\underline{\ad \, P}(-D_\red)^{s,\t}$ of $\ad \, P(-D_\red)^{s,\t}$ consisting of the constant sections, i.e, the kernel of $\nabla \cdot$.  There is thus an exact sequence
\begin{align} \label{nablaconstant}
0 \to \underline{\ad \, P}(-D_\red)^{s,\t} \to \ad \, P(-D_\red)^{s,\t} \xrightarrow{-\nabla \cdot} \ad \, P \otimes K(D_\red)^{s,\t}.
\end{align}
If this were exact at the right, then we would have a resolution of the locally constant sheaf $\underline{\ad \, P}(-D_\red)^{s,\t}$, and so
\begin{align*}
\commsq{ \underline{\ad \, P}(-D_\red)^{s,\t} }{ 0 }{ \ad \, P(-D_\red)^{s,\t} }{ \ad \, P \otimes K(D_\red)^{s,\t} }{ }{ }{ }{ - \nabla \cdot }
\end{align*}
would be a quasi-isomorphism of complexes, in which case we would have
\begin{align*}
\Hy^1(- \nabla \cdot) \cong H^1 \big(X, \underline{\ad \, P}(-D_\red)^{s,\t} \big).
\end{align*}
But $\underline{\ad \, P}(-D_\red)^{s,\t}$ is a local system, so essentially carries only the monodromy information of the connection.  Furthermore, the symplectic form should be recoverable from the cup product
\begin{align*}
H^1 \big( \underline{\ad \, P}(-D_\red)^{s,\t} \big) \otimes H^1 \big( \underline{\ad \, P}(-D_\red)^{s,\t} \big) \to H^2 \big( \underline{\ad \, P}(-D_\red)^{s,\t} \otimes \underline{\ad \, P}(-D_\red)^{s,\t} \big)
\end{align*}
followed by a pairing induced by the Killing form.  However, it does not appear that (\ref{nablaconstant}) is right exact, for $\ad \, P \otimes K(D_\red)^{s,\t}$ contains sections with a $\t$-valued logarithmic term, but any polar term in the image of $-\nabla \cdot$ must be obtained as the bracket with a $\t$-valued form.
\end{rmk}

\section{Families of Bundles Obtained by Hecke Modifications} \label{HM}

The constructions used in the next section depend on an ability to view open sets of the moduli space of bundles as obtained via Hecke modification of a fixed bundle.  We review what this means here.  Further details can be found in \cite{HMWCPB}.

Let $Q$ be a holomorphic principal $G$-bundle over $X$.  A \emph{Hecke modification of $Q$ supported at $x \in X$} consists of a $G$-bundle $P$ and an isomorphism
\begin{align*}
a : P|_{X_0} \arsim Q|_{X_0},
\end{align*}
where $X_0 := X \setminus \{ x \}$.  Let $X_1 \subseteq X$ be an open disc centred at $X$, and choose trivializations of $P$ and $Q$ over $X_1$, and consider the map of trivial bundles over $X_{01} := X_0 \cap X_1$ with respect to these trivializations.  This will yield a holomorphic map $\sigma : X_{01} \to G$.  We can also trivialize $P$ and $Q$ over $X_0$ (any bundle with semisimple structure group is holomorphically trivial over a non-compact Riemann surface); let us do this so that these trivializations correspond via the isomorphism $a$.  Then if $g_{01}$ and $h_{01}$ are the resulting transition functions for $P$ and $Q$, respectively, we have the relationship
\begin{align*}
g_{01} = h_{01} \sigma.
\end{align*}
Essentially, the modification depends only on these maps $\sigma$ up to a suitable equivalence relation.  For $\sigma$ conjugate to a fixed cocharacter $\lambda^\vee$ of $T \subseteq G$, there are finite-dimensional spaces $\Y(Q, \lambda^\vee)$ of equivalence classes of such $\sigma$, which in the absence of more imaginative terminology, will be referred to as spaces of Hecke modifications.  By introducing Hecke modifications at various points in $X$, one obtains families of bundles, with parameter space a (symmetric) product of spaces of Hecke modifications.  The main result of \cite{HMWCPB} can be summed up in the following statement.

\begin{prop} \label{parametrization}
If $G$ is semisimple of adjoint type with root system $A_3, C_l$ or $D_l$ (i.e.\ $G = PGL(4), PSp(2l)$ or  $PSO(2l)$) and if the genus $g$ is even, then one can obtain parametrizations of an open set in the moduli space of principal bundles with spaces of Hecke modifications of the trivial bundle.
\end{prop}

In the case where $G$ is semisimple of adjoint type, if $P$ arises from $Q$ by Hecke modification, then there is a vector bundle $E = E_{PQ}$ of rank $\dim \, G$ and an exact sequence
\begin{align} \label{basseqpb}
0 \to \ad \, P \to E_{PQ} \to E_{PQ}/ \ad \, P \to 0,
\end{align}
so that $E/\ad \, P$ is a torsion sheaf, supported precisely where the Hecke modifications are.  The vector spaces $H^0(X, E/\ad \, P)$ may be identified with the space of infinitesimal deformations of the modification, with the connecting homomorphism
\begin{align*}
H^0(X, E/\ad \, P) \to H^1(X, \ad \, P)
\end{align*}
yielding the Kodaira--Spencer map for the family.  Furthermore, there is also an inclusion of sheaves $E^* \to \ad \, Q$, which is again an isomorphism away from the support of the Hecke modifications, fitting into a commutative diagram
\begin{align} \label{EtoQ}
\vcenter{ \commsq{ E^* }{ \ad \, Q }{ \ad \, P }{ E. }{}{}{}{} }
\end{align}
(See \cite[\S3.2]{HMWCPB}; cf.\ \cite[\S2]{Hurtubise}.)

\section{Isomonodromic Deformation} \label{ID}

\subsection{The Space of Deformation Parameters and the Isomonodromy Connection}

As mentioned in the introduction, isomonodromic deformation has its roots in the Riemann--Hilbert problem, in determining the constraints on the movement of the poles of a connection to ensure the monodromy remains constant.  For $\C\P^1$, these deformation parameters consist of the locations of the poles and the irregular polar part of the connection, but for higher genus the complex structure of the surfaces themselves must also figure in \cite[\S4]{Krichever}; indeed, the complex structures become monodromy parameters themselves.

Since our discussion concerns connections with poles of arbitrary order, we must consider the moduli space of curves with punctures, or marked points, which keeps track of multiplicities.  Recalling our definition of $L$ (\ref{L}), the relevant moduli space will be denoted $\M_{g,L}$ and we now describe it.  We observe that a closed immersion $c : L \to X$ carries the information of both a divisor $D$ (the schematic image of $c$) and a choice of jets (by choosing coordinates at the support of $D$ which pull back precisely to $z \in \C[z]/(z^{k_j})$).  We thus let $\M_{g,L}$ be the space of pairs $(X,c)$, where $X$ is a compact Riemann surface of genus $g$ and $c : L \to X$ is a closed immersion whose image we will take to be a divisor $D$ with a choice of jet of coordinates at each point in its support.

We remark that the space of pairs $(X, D)$, where $D$ is a divisor isomorphic to $L$ but without the information of a choice of jet, is the quotient $\M_{g,L}/ \Aut \, L$.

Recall that to define the monodromy data for a triple $(P, s, \nabla)$, we fixed the irregular part of $\nabla$.  Therefore our space of deformation parameters will be $\M_{g,L} \times (\t_L \cap \h_L)^*$, where as before $(\t_L \cap \h_L)^*$ is realized as Laurent polynomials with coefficients in $\t$ with terms of order at most $-2$.

Following \cite[\S5]{Hurtubise}, define
\begin{align*}
\UUU_{g,L,\varepsilon,T} := \{ (X, c, P, s, \nabla) \, : \, (X, c) \in \M_{g,L}, \, (P, s, \nabla) \in \PPP_X(\varepsilon, D= c(L), T) \}.
\end{align*}
This carries a natural projection to the space of deformation parameters
\begin{align*}
\UUU_{g,L,\varepsilon,T} \to \M_{g,L} \times (\t_L \cap \h_L)^*
\end{align*}
given by
\begin{align*}
(X, c, P, t, \nabla) \mapsto \big(X, c, c^* (s \nabla)_\pol \big).
\end{align*}

Since $c$ induces isomorphisms $G_L \to G_D, T_L \to T_D, H_L \to H_D$, etc., $\UUU_{g, L, \varepsilon, T}$ admits an action of $T_L \cap H_L$ and since this is abelian, the coadjoint action is trivial and hence the map above is invariant under this action, which then gives a projection
\begin{align} \label{monodromybundle}
\UUU_{g,L,\varepsilon,T}/ ( T_L \cap H_L) \to \M_{g,L} \times (\t_L \cap \h_L)^*.
\end{align}
The elements of $\UUU_{g, L, \varepsilon, T}$ then consist of tuples $(X, c, P, s, \nabla)$, where $s$ is a generic compatible framing for $\nabla$ over $D_\red$.  Therefore the fibres of this map are precisely the symplectic spaces $\L_X(\varepsilon, D)_\fr^\eta$ of Section \ref{isoprep}.

What remains to define the monodromy data is a choice of tentacles (Section \ref{tentacles}).  One can imagine that once this choice is made, then by moving around in the space of deformation parameters $\M_{g,L} \times (\t_L \cap \h_L)^*$ phenomena such as $y$ encircling a point of $D$ or one of the $y_j$ crossing over an anti-Stokes direction may occur, rendering the monodromy data computed from our initial choice ill-defined; indeed, such occurrences would be equivalent to another choice of tentacles.  However, what we are concerned with is constructing isomonodromic deformations which are local on the base.  For small changes in the deformation parameters, these phenomena will not occur and so for a given element of $\UUU_{g, L, \varepsilon, T} / (T_L \cap H_L)$, we can fix such a choice and this choice will give us well-defined monodromy data in a neighbourhood of that point.  Therefore, for a sufficiently small neighbourhood $U$ in $\M_{g,L} \times (\t_L \cap \h_L)^*$ we get an isomorphism from the bundle (\ref{monodromybundle}) to $U \times \XX_{g,L}$, and therefore isomonodromic deformation gives a local splitting of the bundle.  This is known as the \emph{isomonodromy connection}.

Since the monodromy maps are symplectic (Proposition \ref{monodromysymplectic}), by composing the monodromy map for one fibre of (\ref{monodromybundle}) with the inverse of another, we get a symplectic identification and hence the following statement.

\begin{prop}
The isomonodromy connection is symplectic.
\end{prop}

\subsection{Constructing Hamiltonians}

In this subsection, we follow closely the arguments used in \cite[\S6]{Hurtubise}.  Given $(P, s, \nabla) \in \L(\varepsilon, D)_\fr^\eta$ lying in the bundle (\ref{monodromybundle}), we may consider its (infinitesimal) isomonodromic deformation $(\widetilde{P}, \tilde{t}, \widetilde{\nabla})^{\text{iso}, v}$ with respect to some tangent vector $v$ to the base $\M_{g,L} \times (\t_L \cap \h_L)^*$, which may be considered as an (infinitesimal) splitting of (\ref{monodromybundle}).

We wish to construct a second (infinitesimal) splitting of (\ref{monodromybundle}) as follows.  We may consider $P$ as obtained from a fixed bundle $Q$ by Hecke modifications as in Section \ref{HM}.  Fixing a connection $\nabla_0$ on $Q$ with a single pole $x_0$ away from $D$ and the support $D_0$ of the Hecke modifications and a framing $t_0$ of $Q$ at $x_0$, we may consider the isomonodromic deformation $(\widetilde{Q}, \tilde{t}_0, \widetilde{\nabla}_0)^{\text{iso}, v}$ of $(Q, t_0, \nabla_0)$ in the analogous bundle (\ref{monodromybundle}).  By ``parallel transport'' of the Hecke modifications along $(\widetilde{Q}, \tilde{t}_0, \widetilde{\nabla}_0)^{\text{iso}, v}$, we get a bundle $\widetilde{P}$ and a framing $\tilde{t}$; considering a deformation of the connection $\nabla$ which ``preserves the Hecke modifications'' uniquely determines a deformation $(\widetilde{P}, \tilde{t}, \widetilde{\nabla})$ of $(P, t, \nabla)$.  These ideas will be made precise in what follows.

We thus obtain two deformations of $(P, t, \nabla)$ in $\UUU_{\varepsilon, D, \fr}$ which lie above the tangent vector $v$ to $\M_{g,L} \times (\t_L \cap \h_L)^*$, and hence the difference must be a tangent vector to the fibre $\L(\varepsilon, D)_\fr^\eta$, which we may recall (\ref{torusirregularhamaction}) is symplectic.  The result is that this is the value of a Hamiltonian vector field on $\L(\varepsilon, D)_\fr^\eta$ at $(P, t, \nabla)$.

\subsubsection{Deformations in the Modulus of the Punctured Riemann Surface}

We will first consider a deformation of the modulus of the punctured curve $c : L \to X$.  The tangent space to $\M_{g,L}$ at $(X, c, D)$ is given by $H^1(X, \T(-D))$, where $\T = \T_X$ is the tangent sheaf to $X$.  Therefore a deformation is given by an element $\mu \in H^1(X, \T(-D))$, which we will realize as a cocycle $(\mu_{\alpha \beta})$ with respect to an open covering described below.  We will also assume that the $\mu_{\alpha \beta}$ are supported near $D$.  The deformed Riemann surface $\widetilde{X} = \widetilde{X}^\mu$ is a complex space over $\Spec \, \C[\epsilon]/(\epsilon^2)$ such that there exists a cartesian diagram
\begin{align} \label{deformationsquare}
\vcenter{ \commsq{ X }{ \widetilde{X} }{ \Spec \, \C }{ \Spec \, \C[\epsilon]/(\epsilon^2). }{}{}{}{} }
\end{align}
Thus, $\widetilde{X}$ has the same underlying topological space as $X$, but over a coordinate patch $U_\alpha$, it has structure sheaf $\OO_{U_\alpha} \otimes_\C \C[\epsilon] / (\epsilon^2) = \OO_{U_\alpha} \oplus \epsilon \OO_{U_\alpha}$.  On the overlaps $U_{\alpha \beta} = U_\alpha \cap U_\beta$, the transition functions must also come with automorphisms of $\OO_{U_{\alpha \beta}} \oplus \varepsilon \OO_{U_{\alpha \beta}}$ which are given by $\C$-derivations $\OO_{U_{\alpha \beta}} \to \OO_{U_{\alpha \beta}}$, i.e.\ which in our case are the vector fields $\mu_{\alpha \beta}$.  So if $\omega$ is a holomorphic $k$-form\footnote{Since we are on a Riemann surface, we will only have $k = 0$ or $k = 1$.} on $X$ defined over $U_{\alpha \beta}$, then we will identify
\begin{align} \label{formidentification}
\omega(z_\beta) = \omega(z_\alpha) + \epsilon \ \mathcal{L}_{\mu_{\alpha \beta}} \omega( z_\alpha),
\end{align}
where $z_\alpha, z_\beta$ are coordinates on $U_\alpha, U_\beta$, respectively and $\mathcal{L}_{\mu_{\alpha \beta}}$ denotes the operation of taking the Lie derivative in the direction of $\mu_{\alpha \beta}$.  In particular, for a function $f$ (i.e.\ a $0$-form),
\begin{align} \label{functionidentification}
f(z_\beta) = f(z_\alpha) + \epsilon df( \mu_{\alpha \beta})(z_\alpha).
\end{align}

An infinitesimal deformation of $(P, t, \nabla)$ is a triple $(\widetilde{P}, \tilde{t}, \widetilde{\nabla})$ over $\widetilde{X}$ whose pullback to $X$ along the upper horizontal map in (\ref{deformationsquare}) is $(P, t, \nabla)$.  A principal bundle $\widetilde{P}$ over $\widetilde{X}$ with a framing $\tilde{t}$ over $D$ is determined by a cocycle in $G( \OO(U_{\alpha \beta}) + \epsilon \OO(U_{\alpha \beta}))$ which is the identity at $D$.  By a connection $\widetilde{\nabla}$ on $\widetilde{P}$, we mean a relative connection over $\C[\epsilon]/(\epsilon^2)$, so given by $1$-forms (over the $U_\alpha$) with values in $\g( \C[\epsilon]/(\epsilon^2))$.  With this, the monodromy data for can be taken and will give a point in $\XX_{g,L}(\C[\epsilon]/(\epsilon^2))$.

The (infinitesimal) isomonodromic deformation (cf.\ \cite[Chapter 0, Definition 16.4]{Sabbah}) of $(P, t, \nabla)$ with respect to the tangent vector $\mu$ will be an infinitesimal deformation $(\widetilde{P}, \tilde{t}, \widetilde{\nabla})^{\text{iso}, \mu}$ of $(P, t, \nabla)$ whose monodromy data lies in the image of
\begin{align*}
\XX_{g,L}(\C) \hookrightarrow \XX_{g,L} \big( \C[\epsilon]/(\epsilon^2) \big).
\end{align*}
Note that if $(\widetilde{P}, \tilde{t}, \widetilde{\nabla})$  is any lift of $(P, t, \nabla)$, then one obtains the monodromy data for $(P, t, \nabla)$ under the map
\begin{align*}
\XX_{g,L} \big( \C[\epsilon]/(\epsilon^2) \big) \hookrightarrow \XX_{g,L}(\C)
\end{align*}
whose composition with the preceding map is an isomorphism.

We now set the notation to explain the details of the constructions described above.  We will let $D = \sum k_j x_j$ be as in (\ref{isoprep}).  Let $D_0 = \{ y_r \}_{r = 1}^s$ be the support of the Hecke modifications and let $x_0$ be the pole of $\nabla_0$.  We will let $U_{1j}$ be a disc centred at $x_j$; $U_{1jl} \subseteq U_{1j}$ will be the Stokes sectors at $x_j$ (the range for the index $l$ will depend on the number of anti-Stokes directions at $x_j$).  Let $U_{2r}$ be a disc centred at $y_r$ and $U_3$ a disc centred at $x_0$.  Assume that the $U_{1j}, U_{2r}, U_3$ are pairwise disjoint.  If $X_1$ is their union, then $X \setminus X_1$ is closed and hence compact and so may be covered by finitely many simply connected open sets $U_{0i}$.  We will use $\alpha, \beta$ for any of these indices.

\paragraph{Isomonodromic Deformation of $(P, t, \nabla)$}

The isomonodromic deformation, $(\widetilde{P}, \tilde{t}, \widetilde{\nabla})^{\text{iso}, \mu}$, of a triple $(P, t, \nabla)$ with respect to $\mu \in H^1(X, \T(-D))$ may be given as follows.  In the open sets $U_{1jl}$ (the Stokes sectors), $U_{2r}, U_3, U_{0i}$---i.e.\ all the open sets which do not intersection $\supp \, D$ and which are by definition simply connected---we may obtain $\nabla$-constant trivializing sections; for the $U_{1jl}$, we may take the fundamental canonical solutions used in defining the Stokes data.  We also choose trivializations over the discs $U_{1j}$ which we may take to agree with the trivializations $t|_{x_j}$ at $x_j$.  Then on all overlaps $U_{\alpha \beta}$ on which we have constant sections, the transition functions $g_{\alpha \beta}$ are constant and all of the monodromy data can be recovered from these functions.  Considering the $g_{\alpha \beta}$ as elements of $G(\OO(U_{\alpha \beta}))$, we may consider their image under the inclusion
\begin{align*}
G\big( \OO(U_{\alpha \beta}) \big) \hookrightarrow G\big( \OO(U_{\alpha \beta}) \oplus \epsilon \OO(U_{\alpha \beta}) \big).
\end{align*}
Clearly, this gives a cocycle which defines a lift to a $G$-bundle $\widetilde{P}$ on $\widetilde{X}$ which pulls back to $P$.  The trivializations on the $U_{1j}$ taken to first order at the $x_j$ give $\tilde{t}$.  Except on the $U_{1j}$, the local connection forms will be zero; on the $U_{1j}$ we can take the same forms to define $\widetilde{\nabla}$.

Instead of $\nabla$-constant trivializations, we may take $\nabla_0$-constant trivializations $t_\alpha$.  These are trivializations of $Q$, but are also trivializations for $P$ away from $D_0$.  Denote the corresponding transition functions by $h_{\alpha \beta}$, which will also be constant.  Where it makes sense, we may write $t_\alpha = s_\alpha \cdot k_\alpha$ for some $G$-valued functions $k_\alpha$; since the $\nabla$-constant and $\nabla_0$-constant trivializations will not coincide, in general, the $k_\alpha$ will not be constant.  With this, we may write
\begin{align*}
g_{\alpha \beta} = k_\alpha h_{\alpha \beta} k_\beta^{-1},
\end{align*}
where the $h_{\alpha \beta}$ are the corresponding transition functions for $Q$.  We will let $A_\alpha := t_\alpha^* \nabla$ be the connection forms for $\nabla$ with respect to these trivializations.  It then follows that
\begin{align*}
A_\alpha = k_\alpha^{-1} dk_\alpha.
\end{align*}
These relations hold over $\widetilde{X}$ as well, however, we should take note of (\ref{functionidentification}), which tells us that
\begin{align*}
k_\beta(z_\beta) =  k_\beta(z_\alpha) + \epsilon dk_\beta \mu_{\alpha \beta} (z_\alpha) = k_\beta( z_\alpha ) \exp( \epsilon \mu_{\alpha \beta} A_\beta(z_\alpha) \big).
\end{align*}
Therefore, the same bundle in the $t_\alpha$ trivializations has transition functions (using the coordinate $z_\alpha$)
\begin{align*}
k_\alpha^{-1} g_{\alpha \beta} k_\beta \exp( \epsilon \mu_{\alpha \beta} A_\beta ) = h_{\alpha \beta} \exp( \epsilon \mu_{\alpha \beta} A_\beta ).
\end{align*}

\paragraph{Isomonodromic Deformation of $(Q, t_0, \nabla_0)$}

The isomonodromic deformation, $(\widetilde{Q}, \tilde{t}_0, \widetilde{\nabla}_0)^{\text{iso}, \mu}$, of $(Q, t_0, \nabla_0)$ can, as above, be described by taking $\nabla_0$-constant trivializations of $Q$ which yield constant transition functions $h_{\alpha \beta}$ and then using those same transition functions to define a bundle $\widetilde{Q}$ over $\widetilde{X}$.

What we want to do now is shift $(P, t, \nabla)$ along, keeping ``constant'' the Hecke modifications we used to construct $P$ from $Q$, i.e.\ obtain a bundle $\widetilde{P}$ which is obtained from $\widetilde{Q}$ by introducing ``the same'' Hecke modifications.  Since the modifications are supported at $D_0$ and we are assuming that the $\mu_{\alpha \beta}$ vanish there, keeping the modifications constant makes sense if we keep the same trivializations (and transition functions).

Now, to obtain the deformation $\widetilde{\nabla}$ of $\nabla$, we will want deformations of the connection matrices $A_\alpha$ that preserve the irregular part of $\nabla$ near $D$ and are holomorphic near $D_0$.  Therefore, we want $\widetilde{\nabla}$ to be represented by matrices
\begin{align*}
A_\alpha + \epsilon a_\alpha
\end{align*}
where the $a_\alpha$ are $\g$-valued 1-forms representing sections of $\ad \, P \otimes K(D_\red)$.  For this to make sense as a connection, we require the usual compatibility condition to hold.  First, recall that (\ref{formidentification}) gives us
\begin{align*}
A_\beta ( z_\beta ) = A_\beta( z_\alpha ) + \epsilon \mathcal{L}_{\mu_{\alpha \beta}} A_\beta( z_\alpha ).
\end{align*}
The compatibility condition is then
\begin{align} \label{deformationcompatibilitycondition} 
A_\alpha(z_\alpha) + \epsilon a_\alpha(z_\alpha) & = \Ad \, h_{\alpha \beta} \big( A_\beta (z_\beta) + \epsilon a_\beta(z_\beta) \big) \nonumber \\
& = \Ad \, h_{\alpha \beta} \bigg( A_\beta(z_\alpha) + \epsilon \big( \mathcal{L}_{\mu_{\alpha \beta}} A_\beta(z_\alpha) + a_\beta(z_\alpha) \big) \bigg). 
\end{align}

We now explain what is meant by a ``deformation of the connection $\nabla$ which preserves the Hecke modifications,'' as alluded to at the beginning of this section.  Dualizing the map $\ad \, P \to E_{PQ}$ of (\ref{basseqpb}) and tensoring with $K(D_\red)$, we obtain the exact sequence
\begin{align} \label{dualexactsequence}
0 \to E^* \otimes K(D_\red) \to \ad \, P \otimes K(D_\red) \to \ad \, P/E^* \otimes K(D_\red) \to 0,
\end{align}
where the last term is a torsion sheaf supported at $D_0$.  Since section of the last term essentially parametrize deformations of the Hecke modifications, for the sections $a_\alpha \in \Gamma(U_\alpha, \ad \, P \otimes K(D_\red))$ to ``preserve'' the modifications, we require them to lie in the kernel of $\ad \, P \otimes K(D_\red) \to \ad \, P / E^* \otimes K(D_\red)$; thus, we will think of them as sections of $E^* \otimes K(D_\red)$.

We now want to prove their existence, i.e., that the $a_\alpha \in \Gamma(U_\alpha, E^* \otimes K(D_\red))$ yielding the $\widetilde{\nabla}$ as suggested above actually exist.  Comparing the coefficient of $\epsilon$ on either side of (\ref{deformationcompatibilitycondition}), this is the same as showing that $\mathcal{L}_{\mu_{\alpha \beta}} A_\beta$ define a coboundary in $\ad \, P \otimes K(D_\red)$.

Observe that since $A_\beta$ is a meromorphic 1-form on a curve, $dA_\beta = 0$ and hence
\begin{align*}
\mathcal{L}_{\mu_{\alpha \beta}} A_\beta = d \iota_{\mu_{\alpha \beta}} A_\beta = d( \mu_{\alpha \beta} A_\beta) = \nabla_0 ( \mu_{\alpha \beta} A_\beta),
\end{align*}
the last since we are using $\nabla_0$-constant frames.

Now, the $A_\beta$ are sections of $\ad \, P \otimes K(D)$ and the $\mu_{\alpha \beta}$ of $\T(-D)$ so the $\mu_{\alpha \beta} A_\beta$ give sections of $\ad \, P$ and $\nabla_0 (\mu_{\alpha \beta} A_\beta)$ sections of $\ad \, P \otimes K$ (or more precisely $\ad \, P \otimes K( m x_0)$, where $m$ is the order of the pole of $\nabla_0$ at $x_0$, but since $\mu_{\alpha \beta} = 0$ near $x_0$, we can ignore this).  We are now thinking of $\nabla_0 ( \mu_{\alpha \beta} A_\beta)$ as a representative of an element of $H^1(X, E^* \otimes K)$ supported near $D$.  If $\zeta \in H^0(X, E(-D_\red)) = H^1(X, E^* \otimes K(D_\red))^*$, we consider its pairing with $\nabla_0( \mu_{\alpha \beta} A_\beta)$:
\begin{align*}
\langle \nabla_0( \mu_{\alpha \beta} A_\beta ), \zeta \rangle = \partial \langle \mu_{\alpha \beta} A_\beta, \zeta \rangle - \langle \mu_{\alpha \beta} A_\beta, \nabla_0 \zeta \rangle.
\end{align*}
Since we are choosing cocycle representatives near $D$, the pairing is given by taking residues at $D$.  But $\langle \mu_{\alpha \beta} A_\beta, \zeta \rangle$ is a (meromorphic) function defined near $D$, so $\partial \langle \mu_{\alpha \beta} A_\beta, \zeta \rangle$ will have no residue.  For the second term, we noted above that $\mu_{\alpha \beta} A_\beta$ has no poles at $D$, and $\nabla_0 \zeta$ will likewise have no residue at $D$.  Therefore
\begin{align*}
\langle \nabla_0( \mu_{\alpha \beta} A_\beta ), \zeta \rangle = 0
\end{align*}
and since $\zeta \in H^0(X, E)$ was arbitrary, it follows that $\mathcal{L}_{\mu_{\alpha \beta}} A_\beta = \nabla_0 ( \mu_{\alpha \beta} A_\beta )$ is a coboundary, and the $a_\alpha$ above exist as claimed.

\

To recap, the isomonodromic deformation of $(P, t, \nabla)$ gave a bundle with transition functions $h_{\alpha \beta} \exp (\epsilon \mu_{\alpha \beta} A_\beta)$ and connection 1-forms $A_\alpha$.  The ``parallel transport'' of the Hecke modifications and irregular part of $\nabla$ yielding $(P, t, \nabla)$ along the isomonodromic deformation of $(Q, t_0, \nabla_0)$ gives the bundle with transition functions $h_{\alpha \beta}$ and connection 1-forms $A_\alpha + \epsilon a_\alpha$ as described above.  Taking the difference gives a $1$-cochain
\begin{align} \label{IMDdiff}
( \mu_{\alpha \beta} A_\beta, -a_\alpha)
\end{align}
which satisfies
\begin{align*}
\nabla (\mu_{\alpha \beta} A_\beta) = d ( \mu_{\alpha \beta} A_\beta ) + [ A_\beta, \mu_{\alpha \beta} A_\beta ] = d ( \mu_{\alpha \beta} A_\beta) = \Ad h_{\alpha \beta}^{-1} a_\alpha - a_\beta,
\end{align*}
noting that $[ A_\beta, \mu_{\alpha \beta} A_\beta] = 0$.  This gives a cocycle in hypercohomology and hence defines a tangent vector to $\L(\varepsilon, D)_\fr^\eta$.

\paragraph{Definition of the Hamiltonian}

We now define a function $H_\mu$ on $\L(\varepsilon, D)_\fr^\eta$ as follows.  We may consider $\nabla_0$ as a connection on $P$ with poles at $D_0$.  Suppose $X_0, X_1$ are as in Section \ref{parametrization} so that if $h_{01}$ is the transition function for $Q$ on $X_{01}$, then $g_{01} = h_{01} \sigma$ is the transition function for $P$, where $\sigma$ is a $G$-valued function on $X_{01}$.  Then if $B_0, B_1$ are the connection matrices for $\nabla_0$, they satisfy
\begin{align*}
B_0 & = \Ad \, h_{01} B_1 - dh_{01} \, h_{01}^{-1} = \Ad \, g_{01} \Ad \, \sigma^{-1} B_1 - dg_{01} g_{01}^{-1} - \Ad \, g_{01} d(\sigma^{-1}) \sigma.
\end{align*}
Therefore, taking
\begin{align*}
A_0 := B_0, \quad A_1 := \Ad \, \sigma^{-1} B_1 - d(\sigma^{-1}) \sigma,
\end{align*}
these define a connection on $P$, which will have poles at $D_0$ with order depending on $\sigma$ and on the root system of $G$.

Therefore, we may think of $\nabla - \nabla_0$ as a section in $H^0(X, \ad \, P \otimes K(D + m D_0))$ for some $m > 0$ and hence
\begin{align*}
\kappa ( \nabla - \nabla_0, \nabla - \nabla_0) \in H^0(X, K^2( 2D + 2m D_0)).
\end{align*}
Since for $(P, s, \nabla) \in \L(\varepsilon, D)_\fr^\eta$, the polar part is fixed, so is that of $\kappa( \nabla - \nabla_0, \nabla - \nabla_0)$ at $D$.  Let $q_0 \in H^0(X, K^2(2D))$ have the same polar part and set
\begin{align*}
q := \kappa(\nabla - \nabla_0, \nabla - \nabla_0) - q_0 \in H^0(X, K^2(D + 2mD_0)).
\end{align*}
Since there is an exact sequence
\begin{align*}
H^1(X, \T(-D - 2mD_0)) \to H^1(X, \T(-D)) \to 0
\end{align*}
we can choose a lift $\tilde{\mu}_{\alpha \beta}$ of $\mu_{\alpha \beta}$ and we define the function on the fibre $\L(\varepsilon, D)_\fr^\eta$ as
\begin{align} \label{hammod}
H_\mu := \tfrac{1}{2} \Res_D \, \tilde{\mu} q,
\end{align}
where again $\Res_D$ means we are summing the residues over $\supp \, D$, and noting that $\tilde{\mu} q \in H^1(X, K)$.  Since we are computing the residue at $D$, $H_\mu$ is independent of the lift $\tilde{\mu}$.

Our aim now is to show that the difference (\ref{IMDdiff}) in the isomonodromic splittings is precisely the Hamiltonian vector field corresponding to $H_\mu$.

Let us consider the open cover $X_0, X_1$ where $X_0 = X \setminus \supp (D + D_0)$ and $X_1$ is a disjoint union of discs centred at the points of $\supp (D + D_0)$.  We write $X_D, X_{D_0}$ for the union of the discs centred at $\supp \, D, \supp \, D_0$, respectively.  Since we have a covering by two open sets, if $(s_{01}, a_0, a_1), (t_{01}, b_0, b_1)$ are two cocycles representing deformations of $(P, t, \nabla)$, then the symplectic form is given by
\begin{align*}
\mbox{\fontsize{10}{12} \selectfont
$\Omega \big( (s_{01}, a_0, a_1), (t_{01}, b_0, b_1) \big) = \Res_{D + D_0} \big( \kappa( s_{01}, \Ad g_{01}^{-1} b_0 + b_1 ) - \kappa( t_{01}, \Ad g_{01}^{-1} a_0 + a_1 ) \big).$ }
\end{align*}
We note that since $X_1$ is a disjoint union of simply connected open sets, we may solve for $\nabla t_{01}' = b_1$ and adjust $t_{01}$ by $t_{01}'$ and thereby assume $b_1 = 0$.

We choose a $\nabla$-constant frame in $X_{D_0}$ and a $\nabla_0$-constant frame in $X_D$.  With this, the pairing of the difference in the isomonodromic deformation (\ref{IMDdiff}) and the deformation $(t_{01}, b_0, 0)$ is
\begin{align*}
\mbox{\fontsize{10}{12} \selectfont
$\Omega \big( (\mu_{01} A_1, a_0, a_1), (t_{01}, b_0, 0) \big) = \Res_{D + D_0} \big( \kappa( \mu_{01} A_1, \Ad g_{01}^{-1} b_0) - \kappa ( t_{01}, \Ad g_{01}^{-1} a_0 + a_1) \big).$ }
\end{align*}
Since $\mu_{01}$ is supported near $D$ and $t_{01}$ is supported near $D_0$, this is
\begin{align*}
\Res_D \, \kappa( \mu_{01} A_1, \Ad g_{01}^{-1} b_0) - \Res_{D_0} \, \kappa ( t_{01}, a_0 + a_1) .
\end{align*}
Since $\Ad g_{01}^{-1} a_0 = a_1 + \nabla ( \mu_{01} A_1)$, the second term is
\begin{align*}
\Res_{D_0} \big( 2\kappa(t_{01}, a_1) + \kappa( t_{01} \nabla (\mu_{01} A_1) ) \big).
\end{align*}
But $a_1$ lies in the image of
\begin{align*}
E^* \otimes K( D_\red) \to \ad \, P \otimes K(D_\red)
\end{align*}
(which is the same as $E^* \otimes K \to \ad \, P \otimes K$ since we are away from $D$), and since $t_{01}$ is represented by sections of $E$, $\kappa(t_{01}, a_1)$ is holomorphic at $D_0$, so the first term vanishes.  The second term vanishes since $\mu_{01}$ is supported at $D$.  Therefore the pairing is given by
\begin{align*}
\Res_D \, \kappa( \mu_{01} A_1, \Ad g_{01}^{-1} b_0).
\end{align*}
But using the fact that in a $\nabla_0$-constant trivialization near $D$, $\nabla - \nabla_0 = A_1$, and $H_\mu = \Res_D \tilde{\mu} q$, and the Leibniz rule, this is precisely the expression for $dH_\mu( t_{01}, b_0, 0)$.

Thus, we have proved the following.

\begin{prop}
For $(P, t, \nabla) \in \L(\varepsilon, D)_\fr^\eta$ and a deformation in $\M_{g,L}$ given by $\mu \in H^1( X, \T(-D))$, the difference between the isomonodromic deformation of $(P, t, \nabla)$ and its parallel transport along that of $(Q, t_0, \nabla_0)$ is the value of the Hamiltonian vector field corresponding to the function $H_\mu$ defined in (\ref{hammod}) on $\L(\varepsilon, D)_\fr^\eta$.
\end{prop}

\begin{rmk} \label{quadraticdifferential}
One will observe that the Hamiltonian defined here is defined in terms of the Killing form, which yields a quadratic Hitchin Hamiltonian, up to a twist by a divisor.  Since quadratic differentials can be thought of as cotangent vectors to the moduli of compact Riemann surfaces, that the Hamiltonian should take such a form makes sense considering that we are lifting a deformation of $\M_{g,L}$.
\end{rmk}

\subsubsection{Deformations in the Irregular Part of the Connection}

In this subsection, we will assume that $Q$ is trivial and $\nabla_0$ is the trivial connection.  This does not limit us in terms of our work so far, since in the cases where parametrizations of the moduli space were obtained in \cite{HMWCPB} (noted in Proposition \ref{parametrization}), $Q$ was taken to be the trivial bundle.

We now consider the isomonodromic deformations induced by a deformation in the irregular part of the connection, i.e.\ in $(\t_L \cap \h_L)^*$.  We will assume that $\supp \, D = \{ x \}$ is a single point and that $\nabla$ has a pole of order $k$ at $x$; accordingly, we will use $D$ to denote the divisor $k \cdot x$.  We can recover the more general situation by considering sums of deformations.  A deformation of the irregular part of the connection is of the form $d\beta$, where
\begin{align} \label{beta}
\beta = \sum_{j= -k+1}^{-1} \beta_j z^j,
\end{align}
is a Laurent polynomial in the coordinate $z$ of order $-k+1$ with values in $\t$ and no holomorphic part.

Here, we will let $\widetilde{X} := X \times_\C \C[\epsilon]/( \epsilon^2 )$; this is the case when $\mu = 0$ in (\ref{deformationsquare}).  The (infinitesimal) isomonodromic deformation of $(P, t, \nabla) \in \L(\varepsilon, D)_\fr^\eta$ with respect to $d\beta$ is an infinitesimal deformation $(\widetilde{P}, \tilde{t}, \widetilde{\nabla})^{\text{iso}, \beta}$ for which the irregular part of $\widetilde{\nabla}$ is $\eta + \epsilon d \beta$.

\paragraph{Isomonodromic Deformation of $(P, t, \nabla)$}

Let $U_1$ be a coordinate disc centred at $x$.  Since we are taking $x$ to be distinct from the pole of $\nabla_0$, we may take a $\nabla_0$-constant trivialization of $Q$ over $U_1$.  Since we are also supposing $x \not\in \supp \, D_0$ (the support of the Hecke modification of $Q$ yielding $P$), $P$ and $Q$ may be assumed to be isomorphic over $U_1$ and hence this $\nabla_0$-constant trivialization also gives one of $P$, which may be taken to agree with the framing $t|_x$.  We will call this the $1$-trivialization, and let $A_1$ be the connection form for $\nabla$ with respect to it.

Over $U_1$ (or passing to a smaller neighbourhood of $x$ if necessary), there exist a holomorphic $G$-valued function $T$ and a meromorphic $\t$-valued $1$-form $B$ such that
\begin{align*}
A_1 = \Ad \, T (B) - dT \cdot T^{-1}.
\end{align*}
Thus, we may think of $T$ as the change of the $1$-trivialization needed to put $\nabla$ in ``diagonal'' form.  We will refer to this trivialization as the $B$-trivialization.  Since Stokes data are computed from the diagonal form, to obtain the isomonodromic deformation, we want to take $\widetilde{\nabla}$ to have connection form
\begin{align*}
B + \epsilon d \beta
\end{align*}
in the $B$-trivialization.  The expression (\ref{beta}) defines $\beta$ as a $\t$-valued function on $U_1$; using the $B$-trivialization, we may consider it as a section $\beta^B$ of $\ad \, P (D - D_\red)$, and in this trivialization
\begin{align*}
(\nabla \beta^B)_B = d \beta + [B, \beta] = d\beta
\end{align*}
since $B, \beta$ commute as they are both $\t$-valued.  Therefore, $\widetilde{\nabla}$ has the expression $(\nabla + \epsilon \nabla \beta^B)_B$ with respect to the $B$-trivialization and hence
\begin{align*}
\widetilde{A}_1 := A_1 + \epsilon (\nabla \beta^B)_1 = A_1 + \epsilon \Ad \, T (d \beta) = A_1 + \epsilon \Ad \, T (\nabla \beta^B)_B
\end{align*}
with respect to the $1$-trivialization.

Now, we may choose the trivializing open cover of $P$ so that on any open set $U_\alpha$ which intersects $U_1$, a $\nabla_0$-constant trivialization exists and hence the resulting transition functions (which are the same for $P$ as for $Q$ since we are away from $D_0$) are constant.  We take $\widetilde{P}$ to be given by $g_{1 \alpha} \exp ( \epsilon (\nabla \beta_\alpha^B))$, where the subscript $\alpha$ means that we are expressing things in terms of the $\alpha$-trivialization.

If $A_\alpha$ is the connection form for $\nabla$ with respect to the $\alpha$-trivialization, then on $U_{1\alpha}$, since the $g_{1\alpha}$ are constant, we have
\begin{align*}
A_1 = \Ad \, g_{1\alpha} \ A_\alpha.
\end{align*}
If we set $\widetilde{A}_\alpha := A_\alpha$, then
\begin{align*}
& \Ad \big( g_{1 \alpha} \exp ( \epsilon (\nabla \beta_\alpha^B)) \big) \widetilde{A}_\alpha - d \big( g_{1 \alpha} \exp ( \epsilon (\nabla \beta_\alpha^B)) \big) \exp ( -\epsilon (\nabla \beta_\alpha^B)) g_{1 \alpha} \\
= & \ \Ad \, g_{1\alpha} \left( A_\alpha + \epsilon \big( [ \beta_\alpha^B, A_\alpha] - d \beta_\alpha^B \big) \right) =  A_1 + \epsilon (\nabla \beta^B)_1,
\end{align*}
which is the compatibility condition for the connection $\widetilde{\nabla}$.  Thus, we get a well-defined infinitesimal deformation $(\widetilde{P}, \tilde{t}, \widetilde{\nabla})$ of $(P, t, \nabla)$.

We want to see that the $(\widetilde{P}, \tilde{t}, \widetilde{\nabla})$ just constructed is in fact the isomonodromic deformation $(\widetilde{P}, \tilde{t}, \widetilde{\nabla})^{\text{iso}, \beta}$ of $(P, t, \nabla)$.  Observe that the connection only has a non-trivial deformation in $U_1$.  Since we can choose the $a$- and $b$-cycles to lie away from $U_1$, the monodromy along these cycles remains unchanged.  Thus, we need only see that the Stokes data at $x$ is also unchanged.  If $U_{1l} \subseteq U_1$ is a Stokes sector, then there exists a $G$-valued function $T_l$ in $U_{1l}$ (the fundamental canonical solution) such that
\begin{align*}
A_1|_{U_{1l}} =: A_{1l} = \Ad \, T_l (B_\pol) - dT_l T_l^{-1}.
\end{align*}
The Stokes factors are given by $T_{l+1}^{-1} T_l$ (up to conjugation by an element of $G$).  But now, if we set $\widetilde{A}_{1l} := A_{1l} + \epsilon (\nabla \beta^B)_{1l}$, then
\begin{align*}
\widetilde{A}_{1l} = \Ad \, T_l(B_\pol + \epsilon d\beta) - dT_l T_l^{-1},
\end{align*}
so it follows that $\widetilde{\nabla}$ has the same fundamental canonical solutions $T_l$ and hence the same Stokes data as $\nabla$.

\paragraph{Isomonodromic Deformation of $(Q, t_0, \nabla_0)$}

Since the deformation $d\beta$ does not change the monodromy data of $(Q, t_0, \nabla_0)$, the isomonodromic deformation $(\widetilde{Q}, \tilde{t}_0, \widetilde{\nabla})^{\text{iso}, \beta}$ is simply the pullback of $(Q, t_0, \nabla_0)$ to $\widetilde{X} = X \times_\C \C[\epsilon]/(\epsilon^2)$ under the natural projection $\widetilde{X} \to X$, so that it has the same transition functions, considered as functions in $G( \OO(U_{\alpha \beta}) \oplus \epsilon \OO(U_{\alpha \beta}))$.  Since $P$ is obtained from $Q$ via Hecke modification, $\widetilde{P}$ will be obtained from $\widetilde{Q}$ by the ``same'' Hecke modifications, and so $\widetilde{P}$ is simply the pullback of $P$ to $\widetilde{X}$.

The deformed connection $\widetilde{\nabla}$ on $\widetilde{P}$ must have local connection forms that satisfy the same compatibility conditions as for $\nabla$ (since $P$ and $\widetilde{P}$ have the same transition functions); such a deformation must therefore be a global one, i.e.\ the local forms must patch together to give a global section of $H^0(X, \ad \, P \otimes K(D))$.

As in (\ref{dualexactsequence}), we have an exact sequence
\begin{align*}
0 \to E^* \otimes K(D) \to \ad \, P \otimes K(D) \to \ad \, P/E^* \otimes K(D) \to 0.
\end{align*}
As just mentioned, a deformation of $\nabla$ will live in $H^0(X, \ad \, P \otimes K(D))$.  For it to also ``preserve the Hecke modifications,'' we will also want it to lie in the image of $H^0(X, E^* \otimes K(D)) \hookrightarrow H^0(X, \ad \, P \otimes K(D))$.


The deformation we want is one which agrees with $\nabla \beta^B$ in a neighbourhood of $x$ (recall that $\nabla \beta^B$ is only defined near $x$).
We will think of $\nabla \beta^B$ as a section in $H^0(X, E^* \otimes K(D)|_D)$ and consider the sequence
\begin{align*}
0 \to E^* \otimes K \to E^* \otimes K(D) \to E^* \otimes K(D)|_D \to 0.
\end{align*}
It comes from a global deformation in $H^0(X, E^* \otimes K(D))$ precisely when its image in $H^1(X, E^* \otimes K) = H^0(X, E)^*$ vanishes.  Using the map $E^* \to \ad \, Q$ from (\ref{EtoQ}), we see that the above sequence in fact fits into a commutative diagram
\begin{align*}
\xymatrix{
0 \ar[r] & E^* \otimes K \ar[r] \ar[d] & E^* \otimes K(D) \ar[r] \ar[d] & E^* \otimes K(D)|_D \ar[r] \ar[d] & 0 \\
0 \ar[r] & \ad \, Q \otimes K \ar[r] & \ad \, Q \otimes K(D) \ar[r] & \ad \, Q \otimes K(D)|_D \ar[r] & 0 }
\end{align*}
whose long exact sequences yield a commutative square
\begin{align*}
\commsq{ H^0(X, E^* \otimes K(D)|_D ) }{ H^1(X, E^* \otimes K) }{ H^0(X, \ad \, Q \otimes K(D)|_D) }{ H^1(X, \ad \, Q \otimes K). }{}{}{}{}
\end{align*}
The first vertical arrow is an isomorphism because $E^* \to \ad \, Q$ is an isomorphism away from the support of the Hecke modification $D_0$, and we are assuming that $D$ and $D_0$ are disjoint.  Furthermore, in \cite[\S4.1]{HMWCPB} it is seen that $H^0(X, \ad \, Q) \to H^0(X, E)$ is an isomorphism, and by Serre duality, the second vertical arrow is as well.  Therefore, the image of $\nabla \beta^B$ in $H^1(X, E^* \otimes K)$ vanishes if and only if its image in $H^1(X, \ad \, Q \otimes K)$ does.  Since we are assuming $Q$, and hence $\ad \, Q$, to be trivial, the map
\begin{align*}
H^0(X, \ad \, Q \otimes K(D)|_D) \to H^1(X, \ad \, Q \otimes K)
\end{align*}
simply takes a meromorphic $\g$-valued differential to its residue.  But since $\nabla \beta^B$ is, by definition, the irregular part of a connection, there is no residue.
It follows that $\nabla \beta^B$ comes from a global section $a \in H^0(X, E^* \otimes K(D))$ as claimed.  We will denote the restriction of $a$ to $U_\alpha$ by $a_\alpha$; by construction, $a_1$ has $\nabla \beta^B$ as its irregular polar part.

The difference between these two isomonodromic deformations is now (in the $\nabla_0$-constant trivialization)
\begin{align*}
( \beta^B, a_0, a_1 - \nabla \beta^B ).
\end{align*}
Let $(t_{01}, b_0, b_1 = 0)$ be an arbitrary deformation as before.  Then the symplectic pairing is given by
\begin{align*}
\Res_{D + D_0} \big( \kappa ( \beta^B, b_0 ) - \kappa ( t_{01}, a_0 + a_1 - \nabla \beta^B ) \big).
\end{align*}
We recall that $t_{01}$ is supported near $D_0$, and $\beta^B$ near $D$, so this simplifies to
\begin{align} \label{hamvfirr}
\Res_D \, \kappa \big( \beta^B, b_0 \big) - \Res_{D_0} \, \kappa( t_{01}, a_0).
\end{align}
But since $a_0$ is a section of $E^* \otimes K(D)$, near $D_0$ we may think of it as a section of $E^* \otimes K$, and $t_{01}$ is represented by sections of $E$, so the second term is the residue of something holomorphic and hence vanishes.

\paragraph{Definition of the Hamiltonian}

Given $(P, t, \nabla)$ and the expression of $\nabla$ as $A$ with respect to a $\nabla_0$-constant trivialization, and its ``diagonalization'' $B$, we define a function $H_\beta$ on $\L(\varepsilon, D)_\fr^\eta$ by
\begin{align} \label{hamirr}
H_\beta := \Res \, \kappa(\beta, B),
\end{align}
the residue being computed at the pole.  We now show that the above deformation (\ref{hamvfirr}) is the value of the Hamiltonian vector field corresponding to this function.

Consider the deformation $(t_{01}, b_0, b_1 = 0)$ above.  The deformation in the connection near $x$ is
\begin{align*}
A + \epsilon b_i .
\end{align*}
Then we must have a corresponding deformation in $B$ satisfying
\begin{align*}
A+ \epsilon b_i & = \Ad \, T \big( \exp(\epsilon t) ( B + \epsilon \tilde{b}_i) \big) = \Ad \, T \big( B + \epsilon ( [t, B] + \tilde{b}_i) \big),
\end{align*}
or
\begin{align*}
\tilde{b}_i = \Ad \, T^{-1} (b_i) - [t, B].
\end{align*}
Therefore, we obtain
\begin{align*}
dH_\beta (t_{01}, b_0, 0) & = \Res \, \kappa( \beta, \Ad \, T^{-1} (b_0) - [t, B] ) = \Res \, \kappa( \Ad \, T (\beta), b_0) + \kappa( [\beta, B], t) \\
& = \Res \, \kappa \big(\Ad \, T (\beta), b_0 \big).
\end{align*}
This proves our claim and the following statement.

\begin{prop}
For $(P, t, \nabla) \in \L(\varepsilon, D)_\fr^\eta$, and a deformation in $(\t_L \cap \h_L)^*$ given by $d\beta \in (\t_L \cap \h_L)^*$, the difference between the isomonodromic deformation of $(P, t, \nabla)$ and its parallel transport along the isomonodromic deformation of $(Q, \nabla_0)$, being the trivial bundle and the trivial connection, is the value of the Hamiltonian vector field corresponding to the function $H_\beta$ defined in (\ref{hamirr}) on $\L(\varepsilon, D)_\fr^\eta$.
\end{prop}

\small

\bibliographystyle{amsalpha}

\end{document}